\documentclass[11pt,draft]{article}
\usepackage[cp1250]{inputenc}
\usepackage{latexsym,amssymb,amsmath}
\title{The commutative Moufang loops with
minimum conditions for subloops I}
\date{}
\author{N. I. Sandu}
\date{}
\begin{document}
\maketitle

\begin{abstract}

The structure of the commutative Moufang loops (CML) with minimum
condition for subloops is examined. In particular it is proved
that such a CML $Q$ is a finite extension of a direct product of a
finite number of the quasicyclic groups, lying in the centre of
the CML $Q$. It is shown that the minimum conditions for subloops
and for normal subloops are equivalent in a CML. Moreover, such
CML also characterized by different conditions of finiteness of
its multiplication groups.
\smallskip\\
 \end{abstract}

The loop $Q$ satisfies \textit{minimum condition for subloops with
the property $\alpha$}, if any decreasing chain of its subloops
with the property $\alpha$ $H_1 \supseteq H_1 \supseteq \ldots$
break, i.e. $H_n = H_{n+1} = \ldots$ for a certain $n$. In this
paper the construction of the commutative Moufang loops
(abbreviated CMLs) with minimum condition for subloops is
examined. In particular, it is shown that such a CML $Q$
decomposes into a direct product of finite number of quasicyclic
groups, which lies in the centre of $Q$, and a finite CML (section
2). In the third section these loops are described with the help
of their multiplication groups. Finally, it is shown in the fourth
section that for the CML, the minimum condition for subloops is
equivalent to the minimum condition for normal subloops, and in
the case of $ZA$-loops these conditions are equivalent to the
minimum condition for normal associative subloops. It follows from
the last statement that the infinite commutative Moufang $ZA$-loop
$Q$ has an infinite centre and if the centre of the CML a
satisfies the minimum condition  for the subloops, then $Q$ itself
satisfies this condition.

We finally denote that loops, in particular the CML, with
different conditions of finiteness are examined in [1 -- 3]. We
remind that the condition of finiteness means such's property,
that holds true for all finite loops, but there exist infinite
loops that do not have this property.

\section{Preliminaries}

Let us bring some notions and results on the theory of the
commutative Moufang loops from [4].
  A \textit{commutative Moufang
loop} (abbreviated CMLs) is characterized by the identity

$$ x^2\cdot yz = xy\cdot xz \eqno{(1.1)}.$$

The \textit{multiplication group} $\frak M(Q)$ of the CML $Q$ is
the group generated by all the \textit{translations} $L(x)$, where
$L(x)y = xy$. The subgroup $I(Q)$ of the group $\frak M(Q)$
generated by all the \textit{inner mappings} $L(x,y) =
L^{-1}(xy)L(x)L(y)$ is called the \textit{inner mapping group} of
the CLM $Q$. The subloop $H$ of the CML $Q$ is called
\textit{normal} (\textit{invariant}) in $Q$, if $I(Q)H = H$.
\smallskip\\
\textbf{Lemma 1.1} [4]. \textit{Let $Q$ be a commutative Moufang
loop with the multiplication group $\frak M$. Then $\frak
M/Z(\frak M)$, where $Z(\frak M)$ is the centre of the group
$\frak M$, and $\frak M^{\prime} = (\frak M,\frak M)$ are locally
finite $3$-groups and will be finite, if $Q$ is finitely
generated}.

The \textit{associator} $(a,b,c)$ of the elements $a, b, c$ of the
CML $Q$ are defined by the equality $ab\cdot c = (a\cdot
bc)(a,b,c)$. The identities

$$L(x,y)z = z(z,y,x) \eqno{(1.2)}$$

$$(x^p,y^r,z^s) = (x,y,z)^{prs} \eqno{(1.3)}$$

$$(x,y,z)^3 = 1, \eqno{(1.4)}$$

$$(xy,u,v) = (x,u,v)((x,u,v),x,y)(y,u,v)((y,u,v),y,x)
\eqno{(1.5)}$$ hold  in the CML [4].

The \textit{centre} $Z(Q)$ of the CML $Q$ is a normal subloop
$Z(Q) = \{x \in Q \vert (x,y,z) = 1 \forall y,z \in Q\}$.
\smallskip\\
\textbf{Lemma 1.2} [4]. \textit{In a commutative Moufang loop $Q$
the following statements hold true: 1) for any $x \in Q$ $x^3 \in
Z(Q)$; 2) the quotient loop  $Q/Z(Q)$ has the index three}.
\smallskip\\
\textbf{Lemma 1.3}  [4]. \textit{The periodic commutative Moufang
loop is locally finite}.
\smallskip\\
\textbf{Lemma 1.4}  [5]. \textit{The periodic commutative Moufang
loop $Q$ decomposes into a direct product of its maximum
$p$-subloops $Q_p$, in addition $Q_p$ belongs to the centre $Z(Q)$
under $p \neq 3$}.

\smallskip

The \textit{system} $\sigma$ of the normal subloops of the loop
$Q$ is called \textit{normal}, if it:

1) contains the loop $Q$ and its identity subloop;

2) is linearly ordered on the inclusion;

3) the intersection and union of any non-empty set of elements of
$\sigma$ is an element of $\sigma$ (fullness). If $A \subseteq B$
are two members of the system $\sigma$ and between them there are
no other members of this system then it is said that the subloops
$A$ and $B$ form a \textit{jump} in the system $\sigma$. The
quotient loop $B/A$ is called the \textit{factor} of this system.
The normal system $\sigma$ is called \textit{central} if for any
jump $A$ and $B$ of the system $\sigma$ $B/A \subseteq Z(B/A)$.
The loop possessing a central system is called a
\textit{$Z$-loop}. This statement is proved in [4, Theorems 4.1,
Chap. VI; 10.1, Chap. VIII].
\smallskip\\
\textbf{Lemma 1.5.} \textit{Any commutative Moufang loop is a
$Z$-loop}.

If the loop possesses a central system entirely ordered by the
inclusion (the \textit{central series}), then this loop is called
\textit{$ZA$-loop}.
\smallskip\\
\textbf{Lemma 1.6} [3]. \textit{Any normal different from the
identity element subloop $H$ of the commutative Moufang $ZA$-loop
$Q$ has a different from identity element intersection with its
centre}.

If the upper central series of the $ZA$-loop have a finite length,
then the loop is called \textit{centrally nilpotent}. The least of
such a length is called the \textit{class} of the central
nilpotentcy.
\smallskip\\
\textbf{Lemma 1.7} [3]. \textit{If a commutative Moufang $ZA$-loop
$Q$ has an infinite associative normal subloop, then its centre
$Z(Q)$ is infinite}.
\smallskip\\
\textbf{Lemma 1.8 (Bruck-Slaby Theorem)} [4]. \textit{The finitely
generated commutative Moufang loop is centrally nilpotent}.
\smallskip\\
\textbf{Lemma 1.9} [3]. \textit{If at least one maximal
associative subloop of the commutative Moufang loop $Q$ satisfies
the minimum condition for subloops, then $Q$ satisfies this
condition itself}.

The CML $Q$ will be called \textit{divisible}, it the equality
$x^n = a$ has at least one solution in $Q$, for any number $n > 0$
and any element $a \in Q$. If $n = 3$, then $a = b^3 \in Z(Q)$ by
Lemma 1.2. Therefore it takes place.
\smallskip\\
\textbf{Lemma 1.10.} \textit{If a subloop of the commutative
Moufang loop $Q$ is divisible, it belongs to the centre $Z(Q)$
and, consequence, is normal in $Q$.}

The \textit{quasicyclic $p$-groups} are some important examples of
divisible CML. As abstract groups they have the set of generators
$1 = a_0, a_1, a_2, \ldots, a_n, \ldots$ and defining relations
$a_0 = a^p_1, a_1 = a^p_2, \ldots, a_n = a^p_{n+1}, \ldots$.

A CML is called \textit{injective} if there exists a homomorphism
$\gamma: B \rightarrow Q$, that $\alpha\gamma = \beta$ for any
monomorphism $\alpha: A \rightarrow B$ and homomorphism $\beta: A
\rightarrow Q$.
\smallskip\\
\textbf{Lemma 1.11.} \textit{The divisible commutative Moufang
loops are injective.}
\smallskip\\
\textbf{Proof.} By Lemma 1.10 a divisible CML is associative, but
divisible abelian groups are injective [6].

Further we will denote by $<M>$ the subloop of loop $Q$, generated
by the set $M \subseteq Q$.
\smallskip\\
\textbf{Proposition 1.12.} \textit{The divisible subloop $D$ of
the commutative Moufang loop $Q$ serves as a direct factor for
$Q$, i.e. $Q = D \times C$ for a certain subloop $C$ of the loop
$Q$. We can choose such a subloop that it possesses the given
before subloop $B$ of the loop $Q$, for which $D \cap B = 1$.}
\smallskip\\
\textbf{Proof.} By Lemma 1.11 there exists such  homomorphism
$\beta: Q \rightarrow Q$, that $\beta\alpha = \varepsilon$ for the
natural inclusion $\alpha: D \rightarrow Q$ and the identity
mapping $\varepsilon: D \rightarrow D$. By Lemma 1.10 the subloop
$D$ is normal in $Q$, therefore $Q = D \times \ker\beta$.

Let now the equality $B \cap D = 1$ hold true for the subloop $B
\subseteq Q$. We denote $H = <D, B>$. By Lemma 1.10 $D \subseteq
Z(Q)$ is the centre of the loop $H$, then it is easy to show that
any element of the CML $H$ has the form $au$, where $a \in B, u
\in D$. By (1.2) and (1.5) we have $L(au,bv)c = c(c,bv,au) =
c(c,b,a) \in B$ for any $a, b \in B$ and any $u, v \in D$.
Consequently, the subloop $B$ is invariant in regard to the inner
mapping group of the CML $H$, hence the subloop $B$ is normal in
$H$. Then $<B, D> = B\times D$ and there is a homomorphism $\xi:
B\times D \rightarrow D$, coinciding with the identity on $D$ and
unitary on $B$. If we replace $\varepsilon$ by $\xi$ in the first
part of this proof, then we obtain $Q = D\times \ker\beta$, where
$B \subseteq \ker\beta$. This completes the proof of Proposition
1.12.

The second part of this proposition states that a divisible CML is
an absolute direct factor.

If the CML $Q$ is given, let us examine the subloop $D$ within it,
generated by all divisible subloops of the CML $Q$. By Lemma 1.10
they all belong to the centre $Z(Q)$ of the CML $Q$, then it is
easy to see that $D$ is a divisible CML. Thus it is the maximal
divisible subloop of the CML $Q$. By the Proposition 1.12 $Q = D
\times C$, where obviously $C$ is a reduced CML, meaning that it
has no non-unitary divisible subloops. Consequently, we obtain
\smallskip\\
\textbf{Proposition 1.13.} \textit{Any commutative Moufanf loop
$Q$ is a direct product of the divisible subloop $D$ that lies in
the centre $Z(Q)$ of the loop $Q$, and the reduced subloop $C$.
The subloop $D$ is unequivocally defined, the subloop $C$ is
defined exactly till the isomorphism.}
\smallskip\\
\textbf{Proof.} Let us prove the last statement. As $D$ is the
maximal divisible subloop of the CML $Q$, it is entirely
characteristic in $C$, i.e. it is invariant in regard to the
endomorphisms of the CML $Q$. Let $Q = D' \times C'$, where $D'$
is a divisible subloop, and $C'$ is a reduced subloop of the CML
$Q$. We denote by $\varphi, \psi$ the endomorphisms $\varphi:Q
\rightarrow D', \psi: Q \rightarrow C'$. As $D$ is an entirely
characteristic subloop $\varphi D$ and $\psi D$ are subloops of
the loop $Q$. It follows from the inclusions $\varphi D \subseteq
D'$ and $\psi D \subseteq C'$ that $\varphi D \cap \psi D = 1$. By
Lemma 1.10 $D$ is a abelian group, therefore $\varphi D, \psi D$
are normal in $D$. Then $d = \varphi d\cdot \psi d$ ($d \in D$)
gives $D = \varphi D \cdot \psi D$, so $D = \varphi D \times \psi
D$. Obviously, $\varphi D \subseteq D \cap D', \psi D \subseteq D
\cap C'$, where from $D = (D \cap D')\times (D \cap C')$. But $D
\cap C' = 1$ as a direct factor of the divisible CML, that is
contained by the reduced CML. Therefore, $D \cap D' \subseteq D, D
\subseteq D'$, i.e. $D = D'$. This completes the proof of
Proposition 1.13.

Let us finally prove.
\smallskip\\
\textbf{Proposition 1.14.} \textit{The following conditions are
equivalent for the commutative Moufang loop $D$:}

\textit{1) $D$ is a divisible loop;}

\textit{2) $D$ is an injective loop;}

\textit{3) $D$ serves as a direct factor for any commutative
Moufang loop that contains it.}
\smallskip\\
\textbf{Proof.} The implication 1) $\longrightarrow 2)$ is proved
in Lemma 1.11.

2) $\longrightarrow$ 3). By the definition of the injective CML
$D$ there is such an homomorphism $\beta: Q \rightarrow D$ that
$\beta\alpha = \epsilon$ for the natural inclusion $\alpha: D
\rightarrow Q$ and identity mapping $\epsilon: D \rightarrow D$.
We denote $\ker \beta = H$. Obviously $Q = <D, H>, H \cap D = 1$
and if $aH = bH$, then $a = b$. Let $x \in Q, d \in D, h \in H$.
The CML is an $IP$-loop, then $(L(x,h)d)H = ((xh)^{-1}(x\cdot hd)H
= (x^{-1}(xd))H = dH$, i.e. $L(x,h)d = d$. Any element from $Q$
has the form $dh$, where $d \in D, h \in H$. Using (1.2) and (1.5)
it is easy to show then that the subloop $D$ is invariant in
regard to the inner mapping group of the CML $Q$, i.e. $D$ is
normal in $Q$. Consequently, $Q = D\times H$.

3) $\longrightarrow$ 1). Let the CML $D$ satisfy the condition 3)
and let there exist such generators $a, b, c$ of the CML $D$, that
$(a,b,c) \neq 1$. Let us examine the CML $Q = <D, x>$, where the
element $x$ does not belong to $D$ and given by all the identity
relations $(a,u,v) = (x,u,v)$ for any $u, v \in D$. Obviously, $D$
is a subloop of the CML $Q$, then it serves as a direct factor.
Therefore the element $x$ associates with any two elements of the
subloop $D$, in particular, $(x,b,c) = 1$. But $(x,b,c) = (a,b,c)
\neq 1$. Contradiction. Consequently, the CML $D$ is associative.
By [6] any abelian group can be embedded as a subgroup into a
divisible group. Therefore the CML $D$ is divisible. This
completes the proof of Proposition 1.14.

\section{Finitely cogenerated commutative  Moufang loops}

A subset $H$ of the CML $Q$ is called \textit{self-conjugate} if
$I(Q)H = H$, where $I(Q)$ is the inner mapping group of the CML
$Q$. A self-conjugate set $L$ of elements of the loop $Q$ will be
called a \textit{normal system of cogenerators}, if any
homomorphism $\varphi: Q \rightarrow H$, for which $L \cap \ker
\varphi \neq$ $ \emptyset$ or $\{1\}$ is a monomorphism, for any
loop $H$. Obviously it is equivalent to the fact that any
non-unitary normal subloop of the loop $Q$ contains an non-unitary
element from $L$.

\smallskip

The loop $Q$ will be called \textit{finitely cogenerated}, if it
possesses a finite normal system of cogenerators.
\smallskip\\
\textbf{Theorem 2.1.} \textit{The following conditions are
equivalent for an arbitrary commutative Moufang loop $Q$:}

\textit{1) $Q$ is a finitely cogenerated loop;}

\textit{2) the loop $Q$ possesses a finite normal subloop $B$,
that $B \cap H \neq \{1\}$, for any normal subloop $H$ of the loop
$Q$;}

\textit{3) the loop $Q$ is a direct product of a finite number of
quasicyclic groups that lie in the centre $Z(Q)$ of the loop $Q$,
and a finite loop;}

\textit{4) the loop $Q$ satisfies the minimum condition for
subloops;}

\textit{5) the loop $Q$ possesses a finite series of normal
subloops any factor of which is either a group of a simple order,
or a quasicyclic group.}
\smallskip\\
\textbf{Proof.} 1) $\longrightarrow 2)$. Let $L$ be a finite
normal system of cogenerators of the CML $Q$ and $a \in Q$ be an
element of an infinite order. By Lemma 1.2 the subloop $<a^{3^n}>$
is normal in the CML $Q$. The intersection $<a^{3^n}> \cap L$ is
either null, or equal to $\{1\}$ for a certain great $n$, that
contradicts the condition 1). Therefore there are no elements of
an infinite order in the CML $Q$. Then, by Lemma 1.3, the subloop
$<L>$ is finite. The system of cogenerators is self-conjugate in
the CML $Q$, then the subloop $<L>$ is normal in $Q$, as the inner
mappings are automorphisms in the CML [4]. Consequently, the
condition 2) holds in the CML $Q$.

2)$\longrightarrow$ 3). It can be shown that the CML $Q$ is
periodic, as it was done when proving the implication 1)
$\longrightarrow$ 2). Then, by Lemma 1.4, it decomposes into a
direct product of its maximal $p$-subgroups, therefore $Q$
contains a finite number of such $p$-subloops. In order to prove
3) we can suppose that $Q, B$ are $3$-loops.

Like in abelian groups [6] the non-negative number $n$, for which
the equality $x^{3^n} = a$ has solutions in $Q$ will be called the
\textit{$3$-height} $h(a)$ of the element $a$. If the equality
$x^{3^n} = a$ has solutions for any $n$, then $a$ will be called
the \textit{infinite $3$-height}, $h(a) = \infty$.

We denote $Q[3] = \{x \in Q \vert x^3 = 1 \}$ and let $a \in
Q[3]$. Then $<\varphi a \vert \varphi \in I(Q)>$ is the minimal
normal subloop containing the element $a$, where $I(Q)$ is the
inner mapping group of the CML $Q$. By the condition 2) $a \in B$,
and then $Q[3]$ will be a finite subloop. It follows from here
that the equality $x^3 = a$ can have not more than a finite number
of solutions in CML $Q$, for an fixed element $a \in Q$. If $h(a)
= \infty$, then the solutions $x_1, \ldots, x_k$ cannot have all a
finite height, as if the equality $y^{3^n} = a$ holds for the
element $y \in Q$, then $y^{3^{n-1}}$ is one of the elements $x_1,
\ldots, x_k$.

Let now $a_1 \in Q[3], h(a_1) = \infty$. We denote the solution of
an infinite height of the equality $a_1 = x^3$ by $a_2$, the
solution of an infinite height of the equality $a_2 = x^3$ by
$a_3$ and so on. Consequently, we have constructed a quasicyclic
group, which lies in the centre of the CML $Q$, by Lemma 1.10,
i.e. it is normal in $Q$. The union $D$ of all quasicyclic groups
of the CML $Q$ is a divisible group, therefore by the Proposition
1.12 $Q = D \times C, D \subseteq Z(Q)$. The subloop $C$ has no
element of an infinite height, as if an element $a \in Q$ of the
order $3^n$ ($n \geq 1$) has an infinite height, then $a^{3^{n-1}}
= a^{3^{-1}}, a^{3^{-1}} \in Q[3]$ and the element $a^{3^{-1}}$
has an infinite height. We have shown that $C[3]$ is a finite
subloop. If $a \in C[3], a^{3^n} = 1, a = x^{3^m}$, then
$x^{3^{n+m}} = 1, x^{3^{m+n-1}} = x^{3^{-1}}, x^{3^{-1}} \in
C[3]$, therefore there is an maximum of heights $k$ of the
elements of subloops $C[3]$. Then $(C[3])^{k+1} = 1$, and by Lemma
1.3 the subloop $C$ is finite. The finiteness of number of the
quasicyclic groups of the CML $Q$  follows from the finiteness of
the subloop $D[3]$.

3) $\longrightarrow$ 4). This statement follows from the fact the
quasicyclic groups and the direct product of their finite number
satisfy the minimum condition for subgroups.

4) $\longrightarrow$ 1). The CML $Q$ has no elements of an
infinite order, as if $a$ is such an element, then $<a^{3^n}>$ ($n
= 1, 2, \ldots$) is a strictly descending series of the subloops
of the CML $Q$. Then, by lemma 1.4, $Q$ decomposes into a direct
product of finite number of maximal $p$-subloops $Q_p$. The
subloop $Q_p[p]$ is normal in $Q$ and it cannot be infinite. In
such a case the subloop $\prod_p Q_p[p]$ will be a finite normal
system of cogenerators.

The implication 3) $\longrightarrow$ 5) follows from Lemma 1.8.

In order to prove the implication 5) $\longrightarrow$ 3) we
should first show that if $Q$ has a finite normal subloop $H$ such
that the quotient loop  $Q/H$ is a quasicyclic group, then $Q$ has
a quasicyclic group of an finite index. First we suppose that the
subloop $H$ is associative. By the definition of the quasicyclic
group of the CML $Q$ is generated by the set $\{a_0H, \ldots,
a_iH, \ldots,\}$, where $a^p_{i+1}H = a_iH, a_0 \in H, i = 1, 2,
\ldots$ We will show that $a_i \in Z_Q(H)$ is the centralizer of
the subloop $H$ in $Q$. If $p = 3$, then if follows from the
equality $a^3_{i+1}h = a_i$, where $h \in H$, for $h_1, h_2 \in H$
from (1.3)- (1.5), that $(a_i,h_1,h_2) = (a^3_{i+1}h,h_1,h_2) =
1$, i.e. $a_i \in Z_Q(H)$. If $p \neq 3$, then by (1.3), (1.4)
$(u^p,v,w) = (u,v,w)^{\pm 1}$. Then we have $(a_1,h_1,h_2) =
(a_1^p,h_1,h_2)^{\pm 1} = (h,h_1,h_2)^{\pm 1} = 1$ from the
relations $a^p_1 = h \in H$. Further, if $a_i \in Z_Q(H)$ and
$a^p_{i+1} = a_ih$, then $(a_{i+1},h_1,h_2) =
(a^p_{i+1},h_1,h_2)^{\pm 1} = (a_ih,h_1,h_2) = 1$ by (1.5), i.e.
$a_{i+1} \in Z_Q(H)$. Therefore $Q = HZ_Q(H)$. As the intersection
$H \cap Z_Q(H)$ is contained in the centre of the CML $Z_Q(H)$,
and the quotient loop $Z_Q(H)/(Z_Q(H) \cap H)$ is isomorphic to
the quasicyclic group $Q/H = Z_Q(H)H/H$ the CML $Z_Q(H)$ is an
infinite abelian group, and it satisfies the minimum condition for
subgroups. Then it contains a quasicyclic group of finite index
[6]. But by the relation $Q = HZ_Q(H)$, the latter has a finite
index in the CML $Q$.

Let now $H$ be an arbitrary subloop. It is finite, then by Lemma
1.8 its upper central series has the form $1 =Z_0 \subset Z_1
\subset \ldots \subset Z_{n-1} \subset Z_n = H$, where
$Z_i/Z_{i-1} = Z(H/Z_{i-1})$ or $Z_i = \{a \in H \vert (a,h_1,
\ldots, h_{2i} = 1 \forall h_1, \ldots, h_{2i}) \in H\}$. (Here
$(u_1,\ldots,u_{2i-1},u_{2i},u_{2i+1}) =
((u_1,\ldots,u_{2i-1}),u_{2i},u_{2i+1})$). The inner map\-pings
are automorphisms in CML [4], then it follows from the last
equality that the subloop $Z_i$ is normal in $Q$, as the subloop
$H$ is normal in $Q$. Further, if follows from the relations

$$Q/H \cong (Q/Z_{n-1})/(H/Z_{n-1}) = (Q/Z_{n-1})/(Z_n/Z_{n-1}) =
(Q/Z_{n-1})/Z(H/Z_{n-1})$$ and according to the previous case that
the CML $Q/Z_{n-1})$ contains a quasicyclic group of finite index.
Without lass of  generality, we will consider that $Q/Z_{n-1}$ is
a quasicyclic group, by the Proposition 1.14. Let us now suppose
that $Q/Z_i$ ($i \leq n - 1$) is a quasicyclic group. Then it
follows from the relations $Q/Z_i \cong (Q/Z_{i-1})/(Z_i/Z_{i-1})
= (Q/Z_{i-1})/Z(Q/Z_{i-1})$ that $Q/Z_{i-1}$ is a quasicyclic
group. We obtain for $i = 1$ that $Q$ contains a quasicyclic group
of finite index.

It is obvious, that the implication 5) $\longrightarrow$ 3) should
be proved supposing that the CML $Q$ contains a series of normal
subloops

$$1 = H_0 \subset H_1 \subset \ldots \subset H_m = Q
\eqno{(2.1)}$$ with $m \geq 2$ that have infinite factors and all
are quasicyclic groups.

Let us show that the series (2.1) contains a member which has a
quasicyclic group of finite index. If the subloop $H_1$ is
infinite, then the statement is obvious. But if it is finite, then
let $H_k$ be such a finite member of the series (2.1) that the
next member $H_{k+1}$ is infinite. Then $H_{k+1}$ contains a
quasicyclic group $L_{k+1}$ of finite index. If all factors of the
series (2.1) which are after the factor $H_{k+1}/H_k$ are finite,
then $L_{k+1}$ has a finite index in $Q$ and by the Proposition
1.14 the statement 3) holds in the CML $Q$.

Let $H_{n+1}/H_n$ be the first infinite factor among  that are
after $H_{k+1}/H_k$. By  Lemma 1.10 the subloop $L_{k+1}$ is
normal in $Q$. There exists a finite normal subloop $H_n/L_{k+1}$
in the CML $H_{n+1}/L_{k+1}$ on which the quotient loop is a
quasicyclic group. By the above proved, the CML $H_{n+1}/L_{k+1}$
contains a quasicyclic group $L_{n+1}/L_{k+1}$ of finite index. In
the CML the quasicyclic groups lie in the centre (Lemma 1.10),
then $L_{n+1}$ is a product of two quasicyclic groups. Continuing
these reasonings, after a finite number of steps we will obviously
obtain that the CML $Q$ contains a subloop that is the direct
product of a finite number of quasicyclic groups of finite index.
Then the CML $Q$ satisfies the condition 3). This completes the
proof of Theorem 2.1.
\smallskip\\
\textbf{Corollary 2.2.} \textit{The commutative Moufang loops
satisfying the minimum condition for subloops, compose a class
closed in regard to the extension.}

The statement follows from the equivalence of the conditions 4)
and 5) of  Theorem 2.1.
\smallskip\\
\textbf{Corollary 2.3.} \textit{The commutative Moufang loops,
satisfying the minimum condition for subloop, are centrally
nilpotent.}

The statement follows from the equivalence of the conditions 3),
4) of the Theorem 2.1 and Lemma 1.8.
\smallskip\\
\textbf{Corollary 2.4.} \textit{The set of elements of any order
is finite in the commutative Moufang loop satisfying the minimum
condition for subloops.}

\section{The multiplication groups of commutative Moufang
loops with minimum condition for subloops}

Let $Q$ be an arbitrary CML and let $H$ be a subset of the set
$Q$. Let $\textbf{M}(H)$ denote a subgroup of the multiplication
group $\frak M(Q)$ of the CML $Q$, generated by the set $\{L(x)
\vert \forall x \in H\}$. Takes place.
\smallskip\\
\textbf{Lemma 3.1.} \textit{Let the commutative Moufang loop $Q$
with the multiplication group $\frak M$, $Z(\frak M)$, which is
the centre of the group $\frak M$ and the centre $Z(Q)$ decompose
into the direct product $Q = D \times H$, moreover, $D \subseteq
Z(Q)$. Then $\frak M = \textbf{M}(D) \times \textbf{M}(H)$,
besides, $\textbf{M}(D) \subseteq Z(\frak M), \textbf{M}(D) \cong
D$.}
\smallskip\\
\textbf{Proof.} It is obvious that any element $a \in Q$ has the
form $a = dh$, where $d \in D, h \in H$. As $d \in Z(Q)$, then
$L(a) = L(d)L(h)$, therefore $\frak M = <\textbf{M}(D),
\textbf{M}(H)>$. It follows from the equality\footnote{A detailed
proof of this equality and the isomorphism $\textbf{M}(D) \cong D$
is presented in author's paper $''$Frattini subloops and
normalizer in commutative Moufang loops$''$.}

$$Z(\frak M) = \{\varphi \in \frak M \vert \varphi = L(a) \forall
a \in Z(Q)\}$$ that $\textbf{M}D) \subseteq Z(\frak M)$, therefore
it is easy to see that the subgroups $\textbf{M}(D),
\textbf{M}(H)$ are normal in $\frak M$ and $\textbf{M}(D) \cong
D$. Finally, if $\varphi \in \textbf{M}(D) \cap \textbf{M}(H)$,
then $\varphi = L(u), L(u)1 \in D \cap H, \varphi$ is an inner
mapping. Consequently, $\frak M = \textbf{D} \times \textbf{H}$,
as required.
\smallskip\\
\textbf{Corollary 3.2.} \textit{The multiplication group $\frak M$
of the periodic commutative Moufang loop $Q$ decomposes into the
direct product of its maximal $p$-subgroups $\frak M_3$, moreover,
$\frak M_p \subseteq Z(\frak M)$ for $p \neq 3$.}
\smallskip\\
\textbf{Proof.} By Lemma 1.4 the CML $Q$ decomposes into a direct
product of its maximal $p$-subgroups, moreover, $Q_p \subseteq
Z(Q)$ for $p \neq 3$. Then it follows from lemma 3.1 that the
group $\frak M$ decomposes into a direct product of the subgroups
$\textbf{M}(Q_p)$, moreover, $\textbf{M}(Q_p) \subseteq Z(\frak
M)$ and $\textbf{M}(Q_p) \cong Q_p$ for $p \neq 3$. In order to
finish the proof, it should be shown that $\textbf{M}(Q_p)$ is a
$3$-group. But this is shown in the next lemma.
\smallskip\\
\textbf{Lemma 3.3.} \textit{The multiplication group $\frak M$ of
the commutative Moufang $3$-loop $Q$ is a $3$-group.}
\smallskip\\
\textbf{Proof.} Let $\gamma$ be an arbitrary element from $\frak
M$. Then $\gamma$ can be presented as a product of a finite number
of translation $\gamma = L(u_1)L(u_2)\ldots L(u_n)$, where $u_1,
u_2, \ldots$ $\ldots, u_n \in Q$. We denote $L = <u_1, u_2,
\ldots, u_n>$. For any element $x \in Q$ we denote by $H(x)$ the
subloop of CML $Q$, generated by set $x \cup L$, by $\frak N(x)$ -
the multiplication group of CML $H(x)$, and by $\Gamma$ - the
subgroup of group $\frak M$, generated by the translations
$L(u_i), i = 1, \ldots, n$. By Lemmas 1.8 and 1.3 $H(x)$ is a
finite centrally nilpotent $3$-loop. Let us show that $\frak N(x)$
is a $3$-loop. Indeed, we denote $H(x) = G$. By Lemma 1.3, Chap.
IV from [4] $\frak M(Z/Z(G)) \cong \frak M(G)/Z^{\ast}$, where
$Z^{\ast} = \{\alpha \in \frak M(G) \vert \alpha x \cdot Z(G) = x
\cdot Z(G) \forall x \in G\}$. If $\theta \in Z^{\ast}$, then we
define the function $f: G \longrightarrow Z(G)$ by the rule
$\theta x = xf(x)$ for $\forall x \in G$. Obviously, $f(x) \in
Z(G)$. If $\eta \in Z^{\ast}$ and $\eta x = xg(x)$, then
$(\theta\eta)x = \theta(L(g(x))x) = L(g(x))\theta x =
(g(x)f(x))x$. Consequently, $Z^{\ast}$ is isomorphic to the group
of one-to-one mappings of CML $Q$ on $Z(G)$. Therefore $Z^{\ast}$
is a $3$-group. If CML $G$ is centrally nilpotent of the class
$k$, then $G/Z(G)$ is centrally nilpotent of class $k - 1$. Then
by inductive assumption $\frak M(G)/Z^{\ast}$ is a $3$-group,
therefore $\frak M(G)$ is also $3$-group.

The restriction $\Gamma$ on $H(x)$ is a homomorphism of $\Gamma$
on the subgroup of the group $\frak N(x)$ which maps the element
$\gamma \in \Gamma$ into the element $L(u_1)\ldots L(u_n)$ from
$\frak N(x)$ of the order $3^t$. Moreover, $\Gamma$ maps $H(x)$
into itself. Consequently, $\gamma^{3^t}$ induces an identity
mapping on $H(x)$. In particular, $\gamma^{3^t}$ maps $x$ into
itself for any $x$ from $Q$. Therefore $\gamma$ has the order
$3^t$. This completes the proof of Lemma 3.3.
\smallskip\\
\textbf{Lemma 3.4.} \textit{The multiplication group $\frak M$ of
an arbitrary commutative Moufang loop is locally nilpotent. But if
group $\frak M$ is periodic, then it is locally finite.}

The proof of the first statement follows from  Lemma 1.1. The
second statement follows from the well-known fact of the group
theory: a periodic locally nilpotent group is locally finite.

Now we can characterize CML, with the minimum condition for
subloops with the help of their multiplication groups.
\smallskip\\
\textbf{Theorem 3.5.} \textit{For an arbitrary non-associative
commutative Moufang loop $Q$ with a multiplication group $\frak M$
the following conditions are equivalent:}

\textit{1) loop $Q$ satisfies the minimum condition for subloops;}

\textit{2) group $\frak M$ is a product of a finite number of
quasicyclic groups lying in the centre of the group $\frak M$, and
a finite group;}

\textit{3) group $\frak M$ satisfies the minimum condition for
subgroup;}

\textit{4) group $\frak M$ satisfies the minimum condition for
normal subgroup;}

\textit{5) group $\frak M$ satisfies the minimum condition for
non-abelian subgroup;}

\textit{6) at least one maximal abelian subgroup of the group
$\frak M$ satisfies the minimum condition for subgroups;}

\textit{7) if group $\frak M$ contain a solvable subgroup of the
class $r$, then $\frak M$ satisfies the minimum condition for
solvable subgroups of the class $r$;}

\textit{8) if group $\frak M$ contain a nilpotent subgroup of the
class $n$, then $\frak M$ satisfies the minimum condition for
nilpotent subgroups of the class $n$.}
\smallskip\\
\textbf{Proof.} 1) $\longrightarrow$ 2). If CML $Q$ satisfies the
minimum condition for subloops, then by Theorem 2.1 $Q = D \times
H$, where $H$ is the direct product of a finite number of
quasicyclic groups, besides, $D \subseteq Z(Q)$, and $H$ is a
finite CML. Then by Lemma 3.1 $\frak M = \textbf{M}(D) \times
\textbf{M}(H)$, and besides $\textbf{D} \subseteq Z(\frak M),
\textbf{M}(D) \cong D$. The group $\textbf{M}(H)$ is finitely
generated, then by Lemma 3.3 $H$ is finite, as it follows from the
Corollary 3.2 that a multiplication group of a periodic CML is
periodic.

The implication 2) $\longrightarrow$ 3) is obvious .Let now the
group $\frak M$ satisfy the condition 3), and the CML $Q$ do not
satisfy the condition 1), and let $Q \supset H_1 \supset H_2
\supset \ldots \supset H_i \supset \ldots$ be an infinite
descending series of subloops of the CML $Q$. It is easy to see
that $\textbf{M}(H_i) \neq \textbf{M}(H_{i+1})$ follows from $H_i
\neq H_{i+1}$, using the relation $\textbf{M}(H_i)1 = H_i$, where
$\textbf{M}(H_i)1 = \{\alpha 1 \vert \alpha \in
\textbf{M}(H_i)\}$. But it contradicts the condition 3).
Consequently, 3) $\longrightarrow$ 1).

By Lemma 3.4 the group $\frak M$ is locally nilpotent, then the
implications 3) $\longleftrightarrow$ 4), 3) $\longleftrightarrow$
5) follow, respectively, from the Theorems 1.24 and Corollary 6.2
from [7].

6) $\longrightarrow$ 3). Let the maximal abelian subgroup $\frak
N$ of the group $\frak M$ satisfy the minimum condition for
subgroups. By Lemma 1.1 the quotient group $\frak M/Z(\frak M)$ is
a $3$-group, therefore by the periodicity of $\frak N$, the group
$\frak M$ is also periodic. Thereof, and in view of the Corollary
3.2, we will consider $\frak M$ a $3$-group. By Lemma 3.4 the
group $\frak M$ is locally nilpotent. Then the condition 6)
$\longrightarrow$ 3) follows from the statement that is proved
using Lemma 1.6, analogous to Theorem 1.19 from [7]:

\textit{it at least one maximal abelian subgroup of the locally
nilpotent $p$-group satisfies the minimum condition for subgroup,
then the group satisfies this condition itself}.

By Lemma 3.4 the group $\frak M$ is locally nilpotent. It is
proved in [8] that for such groups the conditions 3), 7), 8) are
equivalent.

Finally, the implication 3) $\longrightarrow$ 6) is obvious. This
completes the proof of Theorem 3.5.

It is proved in [7] that if the locally finite $p$-group has a
finite maximal elementary abelian subgroup (respect. a finite set
of elements of any order, different from unitary element), then it
satisfy the minimum condition for subgroups (Theorem 1.21
(respect. Theorem 3.2)). Then from the Lemmas 3.3, 3.4 and Theorem
3.6 follows the truth  of the following statement.
\smallskip\\
\textbf{Proposition 3.6.} \textit{The following conditions are
equivalent for an arbitrary commutative Moufang $3$-loop with a
multiplication group $\frak M$:}

\textit{1) the loop $Q$ satisfies the minimum condition for
subloops;}

\textit{2) the group $\frak M$ contains only  a finite set of
elements of a certain order different from the unitary element.}

Finally, let us prove the statement.
\smallskip\\
\textbf{Proposition 3.7.} \textit{The following conditions are
equivalent for an arbitrary non-associative  commutative Moufang
$ZA$-loop $Q$ with a multiplication group $\frak M$:}

\textit{1) the loop $Q$ satisfies the minimum condition for
subloops;}

\textit{2) the group $\frak M$ satisfies the minimum condition for
non-invatiant abelian subgroups.}
\smallskip\\
\textbf{Proof.} Let us first observe that from Lemma 11.4, Chap.
VIII from [4] follows that CML $Q$ is a $ZA$-loop if and only if
its multiplication group is a $ZA$-group.

Let us suppose that the group $\frak M$ satisfies the minimum
condition for non-invatiant abelian subgroups. It follows from the
above-mentioned that it is a $ZA$-group. If $\frak M$ does not
contain non-invatiant abelian subgroups, then, obviously, each
subgroup is normal in it, i.e. it is hamiltonian. However, it is
impossible that the multiplication group of an arbitrary CML
cannot contain a non-abelian  gamiltonian subgroup. Indeed,
arbitrary hamiltonian groups are described by the next theorem
[7]:

An hamiltonian group can be decomposed into a direct product of
the group of quaternions and abelian groups whose each element's
order is not greater that 2. Conversely, a group that has such a
decomposition is hamiltonian.

A \textit{group of quaternions} is the group generated by the
generators $a, b$ and that satisfies the identical relations $a^4
= 1, a^2 = b^2, b^{-1}ab = a^{-1}$. Then it follows from the
Corollary 3.2 that in the case of a multiplication group $a = b =
1$. Consequently, \textit{the arbitrary hamiltonian group of the
multiplication group of CML is abelian.}

Let now $\frak N$ be a non-invatiant abelian subgroup of the group
$\frak M$ and $\alpha$ be an element of infinite order from $\frak
M$. By Lemma 1.1 the quotient group $\frak M/Z(\frak M)$ is a
$3$-group, therefore $\alpha^{3^k} \in Z(\frak M)$ for a certain
natural number $k$. This means that the descending series of
non-invatiant associative subgroups

$$<\frak N, \alpha^{3^k}> \supset <\frak N, \alpha^{3^{k+1}}>
\supset \ldots \supset <\frak N, \alpha^{3^{k+i}}> \supset
\ldots$$ of the group $\frak M$ does not break. But it contradicts
the the condition 5). Consequently, the group $\frak M$ is
periodic. In such a case, we will consider by the Corollary 3.2
that $\frak M$ is a $3$-group.

Let us suppose that the group $\frak M$ does not satisfy the
minimum condition for subgroups. Then, by Lemma 3.4 and Theorem
1.21 from [7] the group $\frak M$ contains an infinite direct
product

$$\frak N = \frak N_1 \times  \frak N_2 \times \ldots \times \frak
N_n \times \ldots$$ of cyclic groups of the order three. If
$\alpha$ is an arbitrary element from the centralizer $Z_{\frak
M}(\frak N)$ of the subgroup $\frak N$ in $\frak M$, then there
exists such a number $n = n(\alpha)$ that

$$<\alpha> \cap (\frak N_{n+1} \times \frak N_{n+2} \times \ldots)
= 1.$$ As the group $\frak M$ satisfies the minimum condition for
non-invatiant abelian subgroups, the infinite descending series of
abelian subgroups

$$\Re^k(\alpha) \supset \Re^{k+1}(\alpha) \supset \ldots ,$$ where
$\Re^k(\alpha) = <\alpha>(\frak N_{k+1} \times \frak N_{k+1}
\ldots )$, contains an non-invatiant subgroup $\Re^k(\alpha)$ ($r
= r(\alpha)$), beginning with a certain natural number $k \geq n$.
As the intersection of all such non-invatiant subgroups coincides
with the subgroup $<\alpha>$, the latter is normal in $\frak M$.
But $\alpha$ is arbitrary element from the centralizer $Z_{\frak
M}(\frak N)$, and it means that $Z_{\frak M}(\frak N)$ is a
hamiltonian group. From here follows that $Z_{\frak M}(\frak N)$
is an abelian group. Obviously, $\frak N_i \subseteq Z_{\frak
M}(\frak N)$, then the minimal subgroup $\frak N_i$ is normal in
$\frak M$. By  Proposition 1.6 from [7], in a $ZA$-group the
minimal normal subgroups are contained in its centre. Then $\frak
N_i \subseteq Z_{\frak M}(\frak N)$, therefore $Z_{\frak M}(\frak
N) = \frak M$. As $Z_{\frak M}(\frak N)$ is an ablian group, the
last equality contradicts the fact that $\frak M$ is an
non-invatiant group. Consequently, the group $\frak M$ satisfies
the minimum condition for subgroups. Then the equivalence of the
conditions 1) and 2) follows from the Theorem 3.5.

\section{The commutative Moufang loops with the minimum
condition for normal subloops}

If it does not cause any misunderstandings, we will further omit
the words ''for subloops'' in the expression ''minimum condition
for subloops''.
\smallskip\\
\textbf{Lemma 4.1.} \textit{Let the series}

$$1 = Z_0 \subset Z_1 \subset \ldots \subset Z_{\alpha} \subset
\ldots \subset Z_{\beta} \subset \ldots \subset Z_{\gamma} = Q
\eqno {(4.1)}$$ \textit{be the upper central series of the
commutative Moufang $ZA$-loop $Q$, $H$ be its arbitrary normal
subloop. Then the non-emptiness of the intersection $H \cap
(Z_{\beta}\backslash Z_{\alpha}$ follows from the non-emptiness of
the intersection $H \cap (Z_{\beta + 1}\backslash Z{_\beta}$ for
any $\beta > \alpha$.}
\smallskip\\
\textbf{Proof.} Let $h \in H \cap (Z_{\beta + 1}\backslash
Z_{\beta}$. The existence of such elements $a, b \in Q$ that
$(h,a,b) \in H \cap (Z_{\beta}\backslash Z_{\alpha})$ follows from
the normality of the subloop $H$ and the definition of the members
of series (4.1). Indeed, if $(h,a,b) \in Z_{\alpha}$ for all $a, b
\in Q$, then $h \in Z_{\alpha + 1} \subset Z_{\beta}$. So, $h
\notin Z_{\beta + 1}\backslash Z_{\beta}$, and it contradicts the
choice of the element $n$. This completes the proof of Lemma 4.1.
\smallskip\\
\textbf{Lemma 4.2.} \textit{Let the commutative Moufang $ZA$-loop
$Q$ be the finite extension of the loop $Q$ satisfies the minimum
condition if and only if the centre $Z(H)$ of the loop $H$ also
satisfies this condition.}
\smallskip\\
\textbf{Proof.} Let us suppose that the centre $Z(Q)$ satisfies
the minimum condition for subloops, and let $a_1, \ldots, a_n$ be
representations of cosets of $Q$ modulo $H$, taken by one from
each coset. We denote $L = <Z(H), a_1, \ldots, a_n>$. Let us show
that the centre $Z(L)$ of the CML $L$ satisfies the minimum
condition. Indeed, the intersection $Z(L) \cap Z(H)$ is contained
into $Z(Q)$, therefore it is a group with minimum condition.
Obviously, the index $Z(H)$ in $<Z(H), Z(L)>$ is finite. We have

$$<Z(H), Z(L)>/Z(H) \cong Z(L)/(Z(L) \cap Z(H)).$$ It follows from
this relation that $Z(L)$ which is a finite extension of the group
$Z(H) \cap Z(L)$, satisfies the minimum condition, satisfies this
condition itself.

Let $N_k$ be a subgroup of the group $Z(L)$ generated by all its
elements, whose order are be divisible by $p^k$. The group $N_k$
is finite, as the group $Z(L)$ satisfies the minimum condition. We
denote by $Z_k$ the subgroup of the group $Z(H)$, generated by all
its elements whose orders are divisors of $p^k$. If $Z(H)$ does
not satisfy the minimum condition, then $Z_k$ should be infinite.
Let

$$Z_k = Z_k^{(1)} \times \ldots \times Z_k^{(m)} \times Z_k^{(m +
1)} \times \ldots $$ be the decomposition of the group $Z_k$ into
an infinite direct product of cyclic groups. If the intersection
$Z_k \cap N_k$ is contained into the finite direct product
$Z_k^{(1)} \times \ldots \times Z_k^{(m)}$, then the intersection
of the groups $M_k = Z_k^{(m + 1)} \ldots$ and $Z(L)$ should
contain only the unitary element:

$$M_k \cap Z(L) = 1. \eqno{(4.2)}$$

By Lemma 1.1 the subgroup $\Phi$ of the inner mapping group of the
CML $Q$ generated by all the mappings of the form $L(a_i,a_J), i,
j = 1, \ldots, n$ is finite. The subgroups $\varphi M_k, \varphi
\in \Phi$ is the set (finite) of all conjugated subloops with
$M_k$ in the CML $Q$, because the elements $a_1, \ldots, a_n$
present a full system of representations of cosets of CML $Q$
modulo $H$, and $M_k \subseteq Z(H)$. The intersection

$$R_k = \cap_{\varphi \in \Phi}\varphi M_k$$ is obviously an
infinite normal subloop in $Q$. We remind that $R_k \subseteq L$,
as $M_k \subseteq Z(L), \varphi M_k \subseteq L$. By Lemma 1.8 and
Lemma 1.6 $R_k \cap Z(L) \neq 1$, that contradicts (4.2).
Consequently, the assumption  that $Z(H)$ does not satisfy the
minimum condition is not  true.

Conversely, let $Z(Q)$ does not satisfy the minimum condition. As
$H$ has a finite index in $Q$, then it follows from the relation

$$Z(Q)H/H \cong Z(Q)/(Z(Q) \cap H)$$ that $Z(Q) \cap H$ has a
finite index in $Z(Q)$. Consequently, $Z(Q) \cap H$ does not
satisfy the minimum condition. But $Z(Q) \cap H \subseteq Z(H)$,
therefore $Z(H)$ does not satisfy the minimum condition  as well.
This completes the proof of Lemma 4.2.
\smallskip\\
\textbf{Lemma 4.3.} \textit{If the commutative Moufang $ZA$-loop,
which is a finite extension of the loop $H$, possesses a normal
subloop $K$, which lies in the centre $Z(H)$ of the loop $H$ and
does not satisfy the minimum condition, then the intersection of
$H$ with the centre $Z(Q)$ of the loop $Q$ does not satisfy the
minimum condition  as well.}
\smallskip\\
\textbf{Proof.} By Lemma 1.4 we'll consider that the CML $Q$ is a
$3$-loop. We denote by $L$ the lower layer of the abelian group
$K$. As $K$ does not satisfy the minimum condition, $L$ is
infinite.

Let us first examine the case when the quotient loop  $Q/H$ is
associative. Let

$$1 = g_1, g_2, \ldots, g_n$$ be a full system of representations
of cosets of CML $Q$ modulo $H$. We suppose by inductive
considerations that the intersection of $L_{i - 1}$ of the centre
of the CML $<H, g_1, \ldots, g_{i - 1}>$ with the subloop $L$ is
infinite. As the quotient loop $Q/H$ is associative, the subloop
$<H, g_1, \ldots, g_{i - 1}>$ inverse image of a normal subloop
under the homomorphism $Q \longrightarrow Q/H$, is normal in $Q$.
The subloop $L$ is invariant in regard to all automorphisms of the
normal subloops $H$ of the CML $Q$. In the CML the inner mappings
are its automorphisms [4]. Then the subloop $L$ is invariant in
regard to the inner mapping group of the CML $Q$, i.e., it is
normal in $Q$. Therefore the intersection $L_{i - 1}$ is also a
normal subloop in $Q$. Let us examine the CML $<L_{i - 1},g_i>$.
By lemma 4.2 this loop's center does not satisfy the minimum
condition. Consequently, if the order of the element $g_i$ is
$3^k$, then there exists such a number $r \leq 3^k$, that for the
infinite set of elements $P$ of the order 3 from the CML $L_{i -
1}$, the elements of the form $pg_i^r, p \in P$ belong to the
centre of the CML $<L_{i - 1}, g_i>$. Now, with the help of (1.1)
we obtain for $p, q \in P$

$$g_i(g_i^rp\cdot g_i^rq) = (g_i\cdot g_i^rp)(g_i^rq),$$

$$g_i(g_i^{2r}\cdot pq) = (g_i^r\cdot g_ip)(g_i^rq).$$

$$g_i^{2r}(g_i\cdot pq) = g_i^{2r}(g_ip\cdot q),$$

$$g_i\cdot pq = g_ip\cdot q.$$ The last equality shows that the
infinite CML $P_i = <P>$ of the index three belongs to the centre
of the CML $<H, g_1, \ldots, g_{i - 1}>$. As $L_{i - 1}$ belongs
to the centre $<H, g_1, \ldots, g_{i - 1}>$ and $P_i \subseteq
P_{i - 1}$, the CML $P_i$ belongs to the centre $<H, g_1, \ldots,
g_{i - 1}, g_i>$. So, the intersection of this CML's centre with
$L_{i - 1}$ is infinite, therefore it does not satisfy the minimum
condition. But $L_{i - 1} \subseteq L_i \subseteq H$, then the
statement is proved in this case.

Let now $Q/H$ be an arbitrary finite CML and by Lemma 1.8  let

$$\overline 1 \subset Z_1/H \subset \ldots \subset Z_k/H =
\overline Q$$ be the upper central series of the CML $Q/H$. By the
first case, the intersection of the centre of the CML $Z_1$ with
the subloop $L$ is infinite. As it has already been proved that
the intersection of the centre of the CML $Z_i$ with the subloop
$L$ is infinite, then applying the first case's results to the CML
$Z_i$ and $Z_{i + 1}$ we obtain that the intersection of the
centre of the CML $Z_{i + 1}$ with the subloop $L$ is also
infinite. For $i + 1 = k$ follows the lemma's statement.
\smallskip\\
\textbf{Lemma 4.4.} \textit{If the periodic commutative Moufang
$ZA$-loop contains an associative normal subloop $H$, that does
not satisfy the minimum condition, then the latter contains a
normal subloop of the loop $Q$, different from itself that does
not satisfy the minimum condition as well.}
\smallskip\\
\textbf{Proof.} By Lemma 1.4 we will consider that $Q$ is a
$3$-loop. We denote by $L$ the lower layer of the group $H$. As
$H$ does not satisfy the minimum condition, $L$ is infinite. Let

$$1 \subset Z_1 \subset \ldots \subset Z_{\gamma} = Q$$ be the
upper central series of the CML $Q$. If  $L \subseteq Z_1$, then
the lemma is proved.

Let us suppose that $L$ does not belong to $Z_1$. The product
$LZ_1 = Q_1$ does not satisfy the minimum condition. The subloop
$L$ is contained in the centre of the CML $Q_1$ ($Q_1 $ is
associative). By the supposition $L$ does not belong to the centre
of the CML $Q$, so, there exists such an ordinal number $\alpha$,
less that $\gamma$ that the centre of the CML $Z_{\alpha}L =
Q_{\alpha + 1}$ does not contain $L$ in its centre anymore.
Consequently, there is such an element $a$ in $Z_{\alpha + 1}$
that the centre $C$ of the finite extension $<Q_{\alpha},a>$ of
the CML $Q$ does not contain the subloop $L$. By Lemma 4.2 the
centre $C$ does not satisfy the minimum condition. The normality
of the subloop $<Z_{\alpha}, a>$ in the CML $Q$ follows from the
relation $Z_{\alpha + 1}/Z_{\alpha} = Z(Q/Z_{\alpha})$, and hereof
follows the normality of the subloop $<Q_{\alpha},a>$.
Consequently, the centre $C$ of the subloop $<Q_{\alpha}, a>$ is
normal in the CML $Q$. By Lemma 4.3 the intersection $C \cap L$
does not satisfy the minimum condition. It is different from the
subloop $L$, as the latter does not belong to $C$. As this
intersection is normal in $Q$, the statement is proved.
\smallskip\\
\textbf{Corollary 4.5.} \textit{In the periodic commutative
Moufang $ZA$-loop $Q$ each associative normal subloop, which
satisfies the minimum condition for the normal subloops of the
loop $Q$ satisfies the minimum condition for its subloops.}

This statement follows from Lemma 4.4.
\smallskip\\
\textbf{Theorem 4.6.} \textit{If at least one maximal associative
subloops of the commutative Moufang $ZA$-loop $Q$ satisfies the
minimum condition for the normal subloops of the loop $Q$, then
$Q$ satisfies the minimum condition for subloops.}
\smallskip\\
\textbf{Proof.} By Lemma 1.2 we will consider that the CML $Q$ is
periodic. Then the statement follows from the Corollary 4.5 and
Lemma 1.9.
\smallskip\\
\textbf{Corollary 4.7.} \textit{In the commutative Moufang
$ZA$-loop the minimum condition for subloops and associative
normal subloops are equivalent.}
\smallskip\\
\textbf{Corollary 4.8.} \textit{If in a commutative Moufang
$ZA$-loop at least one maximal associative normal subloop is
finite, then the loop $Q$ is also finite.}

The statement follows from the Theorems 4.6 and 2.1.
\smallskip\\
\textbf{Corollary 4.9.} \textit{The infinite commutative Moufang
$ZA$-loop $Q$ has an infinite centre.}
\smallskip\\
\textbf{Proof.} By the Corollary 4.8 the CML $Q$ possesses an
infinite associative normal subloop. Then the statement follows
from Lemma 1.7.

We remark that in [4] there is constructed  an example of a CML
with unitary centre.
\smallskip\\
\textbf{Theorem 4.10.} \textit{If the centre $Z(Q)$ of the
commutative Moufang $ZA$-loop $Q$ satisfies the minimum condition
for subloops, then the loop $Q$ satisfies the minimum condition
for subloop itself.}
\smallskip\\
\textbf{Proof.} By the Theorem 2.1 the centre $Z(Q)$ decomposes
into the direct product of a finite number of quasicyclic groups
$D$ and a finite group $C$, and by Proposition 1.12 $Q = D \times
L$. Obviously, the centre $Z(L)$ of the CML $L$ coincides with
$C$. As $C$ is a finite group, then by the Corollary 4.9 the CML
$L$ is finite. Then the CML $Q$ satisfies the minimum condition
for subloops.
\smallskip\\
\textbf{Theorem 4.11.} \textit{If a commutative Moufang loop
satisfy the minimum condition for normal subloops, it satisfies
the minimum condition for subloops as well.}
\smallskip\\
\textbf{Proof.} By Lemma 1.5 an arbitrary CML possesses a central
system. It follows from the minimum condition for normal subloops
that each central system of the CML $Q$ is an ascending central
series, i.e. $Q$ is a $ZA$-loop. Now the statement follows from
the Corollary 4.7.

\smallskip

Tiraspol State University, Moldova

e-mail: sandumn@yahoo.com
\end{document}